\documentclass{gCOM2e}
\usepackage{amsmath,amssymb,graphicx,tikz}

\newcommand{\bfy}{\mathbf y}
\newcommand{\bfb}{\mathbf b}
\newcommand{\tdomain}{\mathcal T}
\newcommand{\ttotal}{$ T_{\text{total}} $}
\newcommand{\tcomm}{$ T_{\text{comm}} $}
\newcommand{\tcpu}{$ T_{\text{cpu}} $}
\newcommand{\taufcpu}{\tau_F^{\text{cpu}}}

\title{A minimal communication approach to parallel time integration}

\author{Andrew T. Barker$^\ast$\thanks{$^\ast$Email: {\tt barker@mpi-magdeburg.mpg.de}.}\\\vspace{6pt}{\em Max Planck Institute for Dynamics of Complex Technical Systems, D-39106 Magdeburg, Germany}}

\begin{document}

\maketitle

\markboth{A.~T. Barker}{International Journal of Computer Mathematics}

\begin{abstract}
We explore an approach due to Nievergelt of decomposing a time-evolution equation along the time dimension and solving it in parallel with as little communication as possible between the processors.  This method computes a map from initial conditions to final conditions locally on slices of the time domain, and then patches these operators together into a global solution using a single communication step.  A basic error analysis is given, and some comparisons are made with other parallel in time methods.  Based on the assumption that parallel computation is cheap but communication is very expensive, it is shown that this method can be competitive for some problems.  We present numerical simulations on graphics chips and on traditional parallel clusters using hundreds of processors for a variety of problems to show the practicality and scalability of the proposed method.
\end{abstract}


\begin{classcode}
65L05, 65M55, 65Y05, 65Y10
\end{classcode}

\section{Introduction}

Consider two different computations.
\begin{equation}
  \textnormal{Use forward Euler to advance a single scalar ODE by a billion timesteps,}
  \label{sequential}
\end{equation}
and
\begin{equation}
  \textnormal{Use forward Euler to advance a billion uncoupled ODEs by one timestep.}
  \label{parallel}
\end{equation}
Implemented in the standard, straightforward way, the problems \eqref{sequential} and \eqref{parallel} require about the same number of floating point operations.  However, each operation in \eqref{sequential} depends on the one before it, while the operations of \eqref{parallel} can all be performed independently.  In particular, \eqref{sequential} relies on the clock speed of a single processor, while \eqref{parallel} depends on the number of processors.  Hardware is moving in the direction of increased concurrency and parallelism, and computational science is moving with it, so that processing units are proliferating and becoming very cheap, while clock speed will remain expensive.  

The catch in the description of \eqref{parallel} is the word ``uncoupled.''  For useful problems in scientific computing, problems are rarely completely uncoupled, which means that different parallel processes will have to communicate somehow.  That leads us to consider a third class of problem, namely
\begin{equation}
  \textnormal{Add a billion numbers, each in distinct (and distant) physical locations.}
  \label{reduce}
\end{equation}
In addition to the floating point operations, problems of type \eqref{reduce} require a great deal of communication, which is already expensive on current hardware.  On future hardware we anticipate that this cost will be the dominant one, and that good algorithms will avoid problems of this type as much as possible.

In this paper we take these trends in hardware seriously, perhaps even taking them to an extreme.  In particular, we assume that computations of type \eqref{parallel} are not just very cheap but essentially free.  We assume that problems of type \eqref{sequential} are more expensive.  Most importantly, we assume that problems of type \eqref{reduce}, namely communication, are {\em very} expensive, so expensive that we will do a thousand or a million problems of type \eqref{parallel} in order to avoid doing one problem of type \eqref{reduce}.

In 1964 Nievergelt presented an approach to the solution of time-dependent differential equations that divides the time domain into slices which are assigned to different processors \cite{nievergeltparallel}.  His approach never attracted much attention because of its high computational cost in terms of floating point operations, but we revisit it here with an eye toward communication cost instead, a factor which was not predicted to be important when Nievergelt was writing.  We present three numerical examples for this method (which we believe are the first actual parallel implementations of the method), and we argue that the Nievergelt approach minimizes total communication among all parallel in time algorithms.

The basic idea of the Nievergelt approach is to {\em construct} a propagation operator on each slice, that is, a map from initial conditions to final conditions, in a way that is embarrassingly parallel but requires a great deal of computation.  These maps are composed with a {\em single} reduce operation of type \eqref{reduce} at the very end of the computation.  As presented here this method has only limited applicability and in particular is difficult to apply to large systems of nonlinear equations.  However, because of its simplicity and its lack of communication, we present it as an instructive example of what kinds of algorithms will minimize communication in future computing environments where parallel computations will become much cheaper in comparison to communication.

\section{Some context}
\label{context}

Most parallel algorithms for evolution equations exploit parallelism in space, often with very impressive speedup and scalability.  Here we focus on parallelism in time, and in particular we will focus on applications where the computational resources are abundant and total time to solution is the primary concern.

One way to exploit parallel computation in the time integration is to use a predictor--corrector framework where the correctors run in parallel one or more steps behind the predictors \cite{christliebparallel}.  The primary purpose of parallelism in this case is to improve the accuracy of the time integration rather than to speed it up, and this method is not suitable for use on a large number of processors.

Some methods based on the matrix exponential are well-suited to parallelism and require minimal communication between processors.  One early effort in this direction is from Gallopoulos and Saad \cite{gallopoulosefficient}, who describe a parallel Krylov subspace method for approximating the matrix exponential.  A more recent approach is the paraexp method of Gander and G\"uttel for linear initial value problems with constant coefficients \cite{ganderparaexp}, which uses a time-domain decomposition together with a ``near-optimal exponential propagator'' to parallelize the solve, with results that show good speedup on up to eight processors.  A similar method is due to Hochbruck and Ostermann \cite{hochbruckexponentialrk}, where again a linear initial value problem can be split into several independent problems, where the division is based on quadrature of the variation-of-constants formula.  These approaches depend on the matrix exponential, so that they are not suitable for problems with a coefficient that varies in time, while we will see that the Nievergelt approach is suitable for such a problem.  Similarly, the contour integral quadrature approach of Sheen et. al. \cite{sheenparallel} relies on the Laplace transform and so cannot be used with varying coefficients, and in addition this method is restricted to parabolic problems.  The similar Dunford-Cauchy integral approach to approximating the matrix exponential in \cite{gavrilyukexponentially} has similar restrictions.  We also mention that none of the above papers include numerical results on more than eight parallel processors.

The parareal scheme has been applied on a large number of processors, is suitable for non-constant coefficient problems, and has been used very effectively to get speedup for a variety of time-dependent problems \cite{lionsresolution,madayparareal,ganderanalysis}.  In this approach the time domain is split into subintervals which are each assigned to processors.  An iteration involving a fine time-stepping scheme (which determines the desired accuracy) and a coarse scheme (which facilitates fast transfer of information between subdomains) is repeated until the error is small.  The parareal methodology has been applied to many types of problems with great success, including for problems in structural and fluid dynamics in \cite{farhattime} and \cite{fischerparareal}, for an ocean model in \cite{liumodified}, and for reaction--diffusion problems in \cite{duarteparareal}, among many other practical and more theoretical studies.  The computational cost of the parareal method can be impressively reduced by combining it with the Spectral Deferred Correction method, but this does nothing to reduce the cost of communication \cite{minionhybrid}.

The approach of Nievergelt that we discuss here is similar to the parareal method in terms of the target applications, namely small systems of differential equations integrated over a very long time interval where spatial parallelization is not easy, but computational resources are abundant.  However, the Nievergelt approach requires the same amount of communication as a {\em single iteration} of the parareal method, measured in terms of total amount of data communicated (or in terms of number of sends and receives).  Although minimal communication is our main goal, we also note that this method will not require a coarse propagator, something that is not always easy to construct.  Finally, the parareal algorithm as well as many of the matrix exponential or Laplace transform based approaches work relatively poorly or not at all for hyperbolic problems \cite{farhattime,ganderanalysismodified,gandercharacteristics}, while we will see that the Nievergelt approach behaves well for the wave equation.

The price we pay for all of these good properties is that the computational cost of Nievergelt is very high when measured simply in terms of floating point operations.  However, almost all of the computations can be done in parallel, so given our (admittedly extreme) assumptions on the computational environment, the Nievergelt approach is competitive for some problems.


In the next section we approach questions of accuracy and scalability in the context of a very simple numerical model problem in order to clarify the properties of the Nievergelt method.  Then in section \ref{sec:costmodel} we present some results on graphics processing units, for which the Nievergelt algorithm is well suited, and discuss the costs and potential parallel speedup of the method in more depth.   In Section \ref{analysis}, we take a slightly more theoretical look at the method and provide some error analysis.  We apply the method to some larger and slightly more realistic problems on standard parallel hardware in Section \ref{pdesection}, and we end with some concluding remarks.

\section{A simple example}
\label{methods}

To make the discussion concrete, consider the model nonlinear initial value problem
\begin{equation}
  y' = y^2, \quad y(0) = 1,
  \label{scalarmodel}
\end{equation}
which has exact solution
\[
  y = \frac{1}{1 - t}.
\]
We will avoid the singularity at time $ t = 1 $ and consider this problem on the time interval $ \tdomain = [0, T] $ with final time $ T = 1/2 $.  We assume the quantity of interest is simply $ y(T) $ and we calculate all our (absolute) errors with respect to that.  

Now we subdivide $ \tdomain $ into subparts $ \tdomain_j = [T_{j-1}, T_j], j = 1 \ldots N $ where $ T_0 = 0, T_N = T $.  On each part $ \tdomain_j $ we consider an initial value problem
\begin{equation}
  y' = y^2, \quad y(T_{j-1}) = \lambda_{j-1}
  \label{localivp}
\end{equation}
on $ \tdomain_j $.  Of course in general $ \lambda_{j-1} $ is unknown, so we conceptually solve the problem for all possible $ \lambda_{j-1} $, constructing not just a scalar value $ y(T_n) $ but a mapping $ \phi^j $ that takes initial values to final values.

We assume $ \lambda_{j-1} $ is known to lie in some range $ \Xi = [a, b] $.  We call $ \Xi $ the {\em initial value space} and sample points inside $ \Xi $ to get a discrete set $ \Xi_h $.  Then, for each $ \xi \in \Xi_h $ we run the initial value problem \eqref{localivp}.  When the initial value $ \lambda_{j-1} $ becomes known (communicated to us from the previous processor in the final step of the algorithm), we use interpolation to calculate $ \lambda_j $.

\begin{table}
  \tbl{Errors for the Nievergelt algorithm on four time domains, with different discretization sizes in time $ \Delta t $ and varying number $ M $ of Chebyshev interpolation in the initial value space, using backward Euler.}
  {\begin{tabular}{l|lllll}
      & $ M $& & & & \\
          $ \Delta t $&         3&       4&        5&         6&   7 \\
      \colrule
      0.01&   0.0681& 0.0201& 0.0274& 0.0278& 0.0278 \\ 
      0.005&  0.0511& 0.00751& 0.0138& 0.0141& 0.0141 \\ 
      0.0025& 0.0422& 0.00098& 0.00672& 0.00700& 0.00698 \\ 
      0.001&  0.0370& 0.00293& 0.00254& 0.00280& 0.00278 \\ 
      0.0001& 0.0339& 0.00526& 0.000050& 0.000296& 0.000278 \\ 
    \end{tabular}}
  \label{foursubdomains}
\end{table}

\begin{figure}
  \begin{center}
   \includegraphics[width=0.45\textwidth]{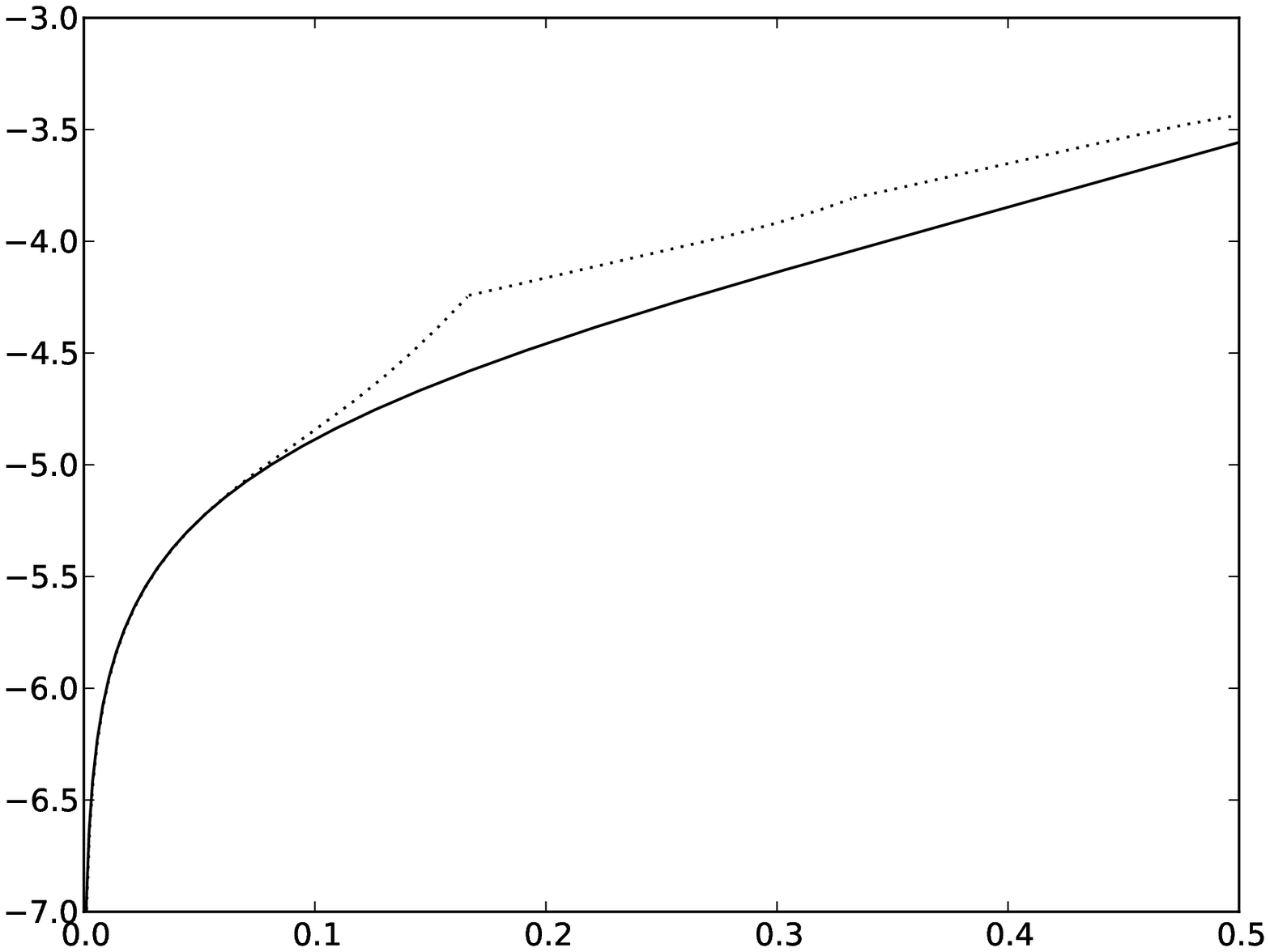}
   \includegraphics[width=0.45\textwidth]{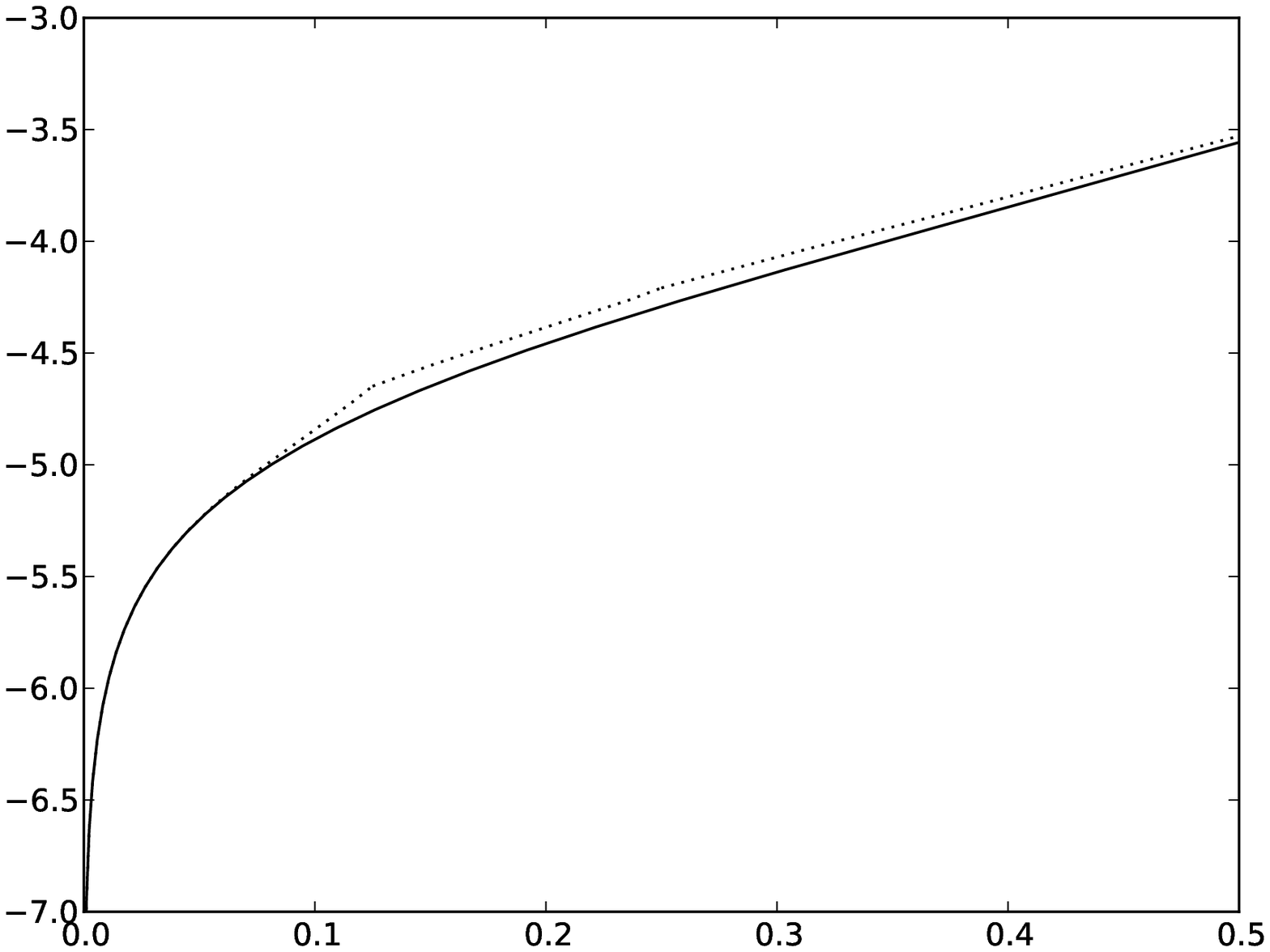}
    \caption{The error of the Nievergelt approach with $ \Delta t = 10^{-4} $ on a log scale on the vertical axis plotted against time on the horizontal axis.  The solid line is the error with plain backward Euler timestepping and no parallel framework, while the dotted line shows the error using $ M = 6 $ Chebyshev points.  The left figure has 3 parallel time intervals, and the right has 4 parallel time intervals.  The interpolation does add to the error, but the amount is comparable to the original time discretization error.}
    \label{errorpartition}
  \end{center}
\end{figure}

To understand how the interpolation process affects accuracy, we represent the initial value space $ \Xi = [0,2] $ as a Lagrange interpolant on the $ M $ Chebyshev nodes in this interval.  We present some accuracy results in Table \ref{foursubdomains}.  Intuitively we wish to balance the time--discretization errors with interpolation errors that are of the same order to get an efficient method.  Since Chebyshev interpolation is spectrally accurate, we see in Table \ref{foursubdomains} that we need just 5 or 6 Chebyshev points is enough to get $ 10^{-4} $ accuracy.  Also see Figure \ref{errorpartition}, which shows the component of error related to interpolation as the time integration proceeds.  We remark that the number of subdomains $ N $ does not significantly affect these results.

\begin{table}
  \tbl{Errors for the Nievergelt approach and the parareal approach, with $ \Delta t = 0.0001 $ and various numbers of processors $ N $.  For Nievergelt we use $ M = 6 $ Chebyshev nodes for interpolation.  Both approaches use backward Euler for the fine time discretization, and the parareal results also use backward Euler with $ \Delta T = 0.1 $ for the coarse propagator.}
  {\begin{tabular}{l|llll}
     &        & parareal& parareal& parareal \\
  $N$& Nievergelt& $k=2$&  $k=3$&  $k=5$ \\
  \colrule
  1&   2.77e-4& 2.77e-4& 2.77e-4& 2.77e-4 \\ 
  2&   8.50e-4& 2.77e-4& 2.77e-4& 2.77e-4 \\
  4&   2.83e-4& 9.46e-3& 7.41e-5& 2.77e-4 \\
  8&   2.79e-4& 1.27e-2& 2.50e-4& 2.77e-4 \\
  16&  2.72e-4& 3.18e-3& 2.13e-4& 2.77e-4 \\
  32&  2.54e-4& 9.49e-4& 2.68e-4& 2.76e-4 \\
  64&  2.16e-4& 4.36e-4& 2.73e-4& 2.74e-4 \\
  \end{tabular}}
  \label{errorcomparison}
\end{table}

\begin{figure}
  \begin{center}
\begin{tikzpicture}[scale=1.0]
    \draw[->] (-0.5, 0.0) -- node[below=0.5cm] {processor} (6.5, 0.0);
    \draw[->] (0.0, -0.5) -- node[left] {time} (0.0, 4.5);

    \foreach \x/\xtext in {0.0/0, 1.0/1, 2.0/2, 3.0/3, 4.0/4, 5.0/5}
      {\filldraw[fill=lightgray,draw=black,xshift=\x cm] (0.0, 0.0) rectangle (1.0, 3.0);
      \draw[xshift=0.5cm] (\x, 0.0) node[below] {\xtext};}

    \foreach \x/\y in {1.0/3.0, 2.0/3.25, 3.0/3.5, 4.0/3.75, 5.0/4.0}
      {\filldraw[black,xshift=\x cm,yshift=\y cm] (0.0, 0.0) rectangle (1.0, 0.25);
      \filldraw[fill=white,draw=black] (\x,\y) circle (0.1);}
\end{tikzpicture}
  \end{center}
  \caption{Timing and communication in the Nievergelt approach.  In this figure the gray squares represent construction of the propagation operator (usually done by timestepping across $ \tdomain_j $ with several different initial conditions), the black squares represent computing the interpolation of this operator at an initial condition communicated from the previous processor, and white circles are communication steps.  Compare with Figure \ref{pararealscheduling}.}
  \label{minimalscheduling}
\end{figure}
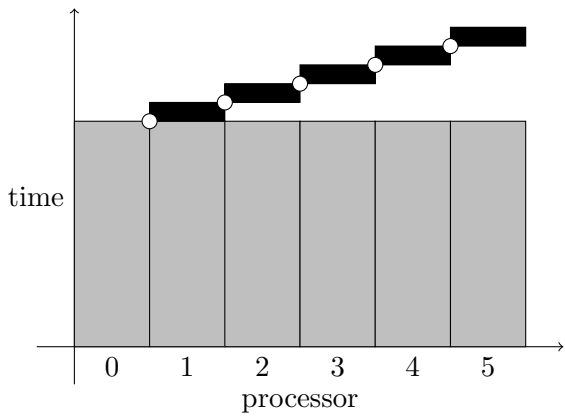

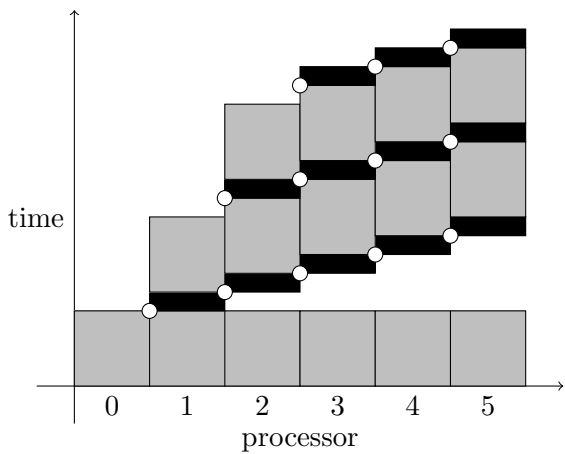
\begin{figure}
  \begin{center}
\begin{tikzpicture}[scale=1.0]
    \draw[->] (-0.5, 0.0) -- node[below=0.5cm] {processor} (6.5, 0.0);
    \draw[->] (0.0, -0.5) -- node[left] {time} (0.0, 5.0);

    \foreach \x/\xtext in {0.0/0, 1.0/1, 2.0/2, 3.0/3, 4.0/4, 5.0/5}
      {\filldraw[fill=lightgray,draw=black,xshift=\x cm] (0.0, 0.0) rectangle (1.0, 1.0);
      \draw[xshift=0.5cm] (\x, 0.0) node[below] {\xtext};}

    \foreach \x/\y in {1.0/1.0, 2.0/1.25, 3.0/1.5, 4.0/1.75, 5.0/2.0}
      {\filldraw[black,xshift=\x cm,yshift=\y cm] (0.0, 0.0) rectangle (1.0, 0.25);}

    \foreach \x/\y in {1.0/1.25, 2.0/1.5, 3.0/1.75, 4.0/2.0, 5.0/2.25}
      {\filldraw[fill=lightgray,draw=black,xshift=\x cm,yshift=\y cm] (0.0, 0.0) rectangle (1.0, 1.0);}

    \foreach \x/\y in {2.0/2.5, 3.0/2.75, 4.0/3.0, 5.0/3.25}
      {\filldraw[black,xshift=\x cm,yshift=\y cm] (0.0, 0.0) rectangle (1.0, 0.25);}

    \foreach \x/\y in {2.0/2.75, 3.0/3.0, 4.0/3.25, 5.0/3.5}
      {\filldraw[fill=lightgray,draw=black,xshift=\x cm,yshift=\y cm] (0.0, 0.0) rectangle (1.0, 1.0);}

    \foreach \x/\y in {3.0/4.0, 4.0/4.25, 5.0/4.5}
      {\filldraw[black,xshift=\x cm,yshift=\y cm] (0.0, 0.0) rectangle (1.0, 0.25);
       \filldraw[fill=white,draw=black] (\x,\y) circle (0.1);}

    \foreach \x/\y in {1.0/1.0, 2.0/1.25, 3.0/1.5, 4.0/1.75, 5.0/2.0}
      {\filldraw[fill=white,draw=black] (\x,\y) circle (0.1);}

    \foreach \x/\y in {2.0/2.5, 3.0/2.75, 4.0/3.0, 5.0/3.25}
      {\filldraw[fill=white,draw=black] (\x,\y) circle (0.1);}
\end{tikzpicture}
  \end{center}
  \caption{Timing and communication in the parareal method.  In this figure the gray squares represent action of the fine level timestepping scheme, the black squares represent the action of the coarse timestepping scheme, and white circles are communication steps.  This figure shows $ k = 3 $ iterations of parareal.  Compare with Figure \ref{minimalscheduling}.}
  \label{pararealscheduling}
\end{figure}

For comparison we also solve this problem using a parareal method \cite{madayparareal}.  The parareal algorithm proceeds by iteratively applying the correction equation
\begin{equation}
  \lambda_{j+1}^{k+1} = G_{\Delta T}^j(\lambda_j^{k+1}) + g_{\Delta t}^j(\lambda_j^k) - G_{\Delta T}^j(\lambda_j^k),
  \label{pararealcorrection}
\end{equation}
where $ \lambda_j^k $ is the initial value for time slice (or processor) $ j $ in iteration $ k $, $ g_{\Delta t} $ is a standard time stepping method with a short timestep, and $ G_{\Delta T} $ is a (cheap) coarse time stepping method with larger timestep $ \Delta t $.  To understand the two compared algorithms visually, see Figures \ref{minimalscheduling} and \ref{pararealscheduling}, where we present the scheduling of the Nievergelt and parareal algorithms as in the figures in \cite{minionhybrid}.  Here we have made an analogy between the fine propagator of parareal and the construction of the map from initial to final conditions in the Nievergelt approach, but we should note that this latter operation is much more expensive (by a constant factor of $ M $).  However, like the fine propagation of parareal, the construction of this map is embarrassingly parallel, requiring no synchronization.  Also, we have compared the coarse propagator of parareal to the interpolation in Nievergelt.  These are quite different operations, and it is not obvious exactly how their costs compare, but they play similar roles in the algorithm.  It is clear from these figures that Nievergelt involves significantly fewer communications, especially as the number of parareal iterations increases.

The results of the comparison are given in Table \ref{errorcomparison}.  Here we use backward Euler as both the fine and coarse propagator, with fine discretization size $ \Delta t $ and coarse discretization size $ \Delta T $, and $ k $ represents the number of parareal iterations used.  More iterations gives better accuracy, but recall that each iteration of parareal requires as much communication as the entire Nievergelt algorithm, and in Table \ref{errorcomparison} we see that even for this simple problem it usually takes a few iterations to get to discretization error.  

To end this section we summarize the Nievergelt approach in pseudocode in Figure \ref{pseudocode}.

\begin{figure}
\begin{center}
\fbox{\begin{minipage}{0.8\textwidth}
\begin{raggedright}
  {\bf Given.} An initial value problem $ y' = f(t,y(t)), y(0) = y_0 $, a sample $ \Xi_h $ of the initial value space, a final time $ T $, and $ N $ processors. \\
  {\bf Initialize.} Let $ T_j = jT/N $ and assign the interval $ [T_{j-1}, T_j] $ to processor $ j $.\\
  {\bf Computation phase.} \\
  \hspace*{.5cm} {\bf for each} processor $ j $: \\
  \hspace*{1cm} {\bf for each} $ \xi_k \in \Xi_h $: \\
  \hspace*{1.5cm} solve $ y' = f(t,y(t)), y(T_{j-1}) = \xi_k $ \\
  \hspace*{1.5cm} let $ \lambda_j^k = y(T_j) $ \\
  \hspace*{1cm} construct approximating function $ \phi(\xi) $ such that $ \phi(\xi_k) = \lambda_j^k $. \\
  {\bf Communication phase.} \\
  \hspace*{.5cm} {\bf for each} processor $ j $: \\
  \hspace*{1cm} receive $ \lambda_{j-1} $ from previous processor \\
  \hspace*{1cm} calculate $ \lambda_j = \phi(\lambda_{j-1}) $ \\
  \hspace*{1cm} send $ \lambda_j $ to the next processor\\
\end{raggedright}
\end{minipage}}
\end{center}
\caption{The minimal communication Nievergelt approach in pseudocode.}
\label{pseudocode}
\end{figure}

\section{A graphics card implementation and cost model}
\label{sec:costmodel}

\begin{table}
  \tbl{Time to solution in microseconds and speedup of the GPU implementation of the Nievergelt method as compared to a simple serial implementation on the CPU.}
{\begin{tabular}{lll|ll|l}
  $\Delta t$&  $N$&  $M$&  \ttotal& \tcpu / \ttotal& $T_{\text{estimate}}$ \\
  \colrule
  $2^{-14}$&  32&   4&   193.6&  2.2&  200.4 \\
  $2^{-14}$&  64&   4&   187.5&  2.3&  202.5 \\
  $2^{-14}$&  128&  4&   208.4&  2.0&  237.3 \\
  $2^{-15}$&  32&   5&   253.7&  3.3&  261.2 \\
  $2^{-15}$&  64&   5&   228.5&  3.8&  233.0 \\
  $2^{-15}$&  128&  5&   257.8&  3.4&  252.5 \\
  $2^{-16}$&  32&   7&   388.3&  4.4&  443.7 \\
  $2^{-16}$&  64&   7&   338.9&  5.0&  324.2 \\
  $2^{-16}$&  128&  7&   392.0&  4.3&  298.1 \\
\end{tabular}}
\label{tab:gpuresults}
\end{table}

One possible use case for the Nievergelt method is on general purpose graphical processing units, which are highly parallel and energy efficient computing units that are widely and cheaply available due to their use in gaming hardware.  Effectively using this hardware is a challenge because they are most efficient when the computational task can be broken down into very many independent threads, which is precisely what the Nievergelt algorithm does.  In Table \ref{tab:gpuresults} we show the speedup of using the parallel capabilities of the GPU over the same computation performed in serial on the CPU, with times recorded in microseconds, for the nonlinear scalar initial value problem, for various timestep sizes $ \Delta t $.  We partition the time interval into $ N $ time slices, for each of which we compute the map from initial to final conditions independently.  For each of the $ N $ time slices we propagate $ M $ initial conditions so that we can use $ MN $ completely independent threads on the GPU.  The computation is finished by doing some cheap interpolation on the CPU.  We are using backward Euler for the timestepping and Chebyshev interpolation with $ M $ nodes (chosen to insure that the total error in the Nievergelt method is comparable to time discretization error) and doing all computations in single precision.  In each case the error for the Nievergelt implementation is comparable to the discretization error due to timestepping.  The timing results include all the time required to allocate memory, build the initial value arrays, transfer data to and from the GPU, and to perform the interpolation that evaluates the map from initial to final conditions.  The CPU and GPU tests are run on the same machine, which has a six core AMD Opteron CPU and an nVidia Tesla C2075 GPU.  This experiment demonstrates the potential for using the inherent parallelism of the graphics card to solve the problem faster than is possible with the CPU of the same machine, even though many more floating point operations are required.


\begin{table}
  \tbl{Approximate measures of parameters in cost models \eqref{gpucostmodel}, \eqref{serialcostmodel}, in microseconds, computed as a least--squares fit from the same computational experiments used to make Table \ref{tab:gpuresults}}{
  \begin{tabular}{lll}
    cost& symbol& value \\
    \hline
    advance one timestep (GPU)& $\tau_F$&    0.040  \\
    apply the Nievergelt map& $\tau_N$&    0.701  \\
    GPU overhead& $\tau_K$&    137.   \\
    advance one timestep (CPU)& $\taufcpu$&  0.051  \\
  \end{tabular}}
  \label{tab:tauparameters}
\end{table}

One possible model for the cost (in terms of time to solution) for the Nievergelt method using the GPU implementation is
\begin{equation}
  C_G = \frac{M n}{N} \tau_F + N \tau_N + \tau_K
  \label{gpucostmodel}
\end{equation}
where $ \tau_F $ is the time to advance a single timestep using backward Euler for a single initial value problem, $ \tau_N $ is a measure of the cost of applying the map that has been constructed (the final communication and interpolation step), and $ \tau_K $ represents some fixed GPU startup costs.  Predicting the time to solution for the experiments in Table \ref{tab:gpuresults} using parameter values from Table \ref{tab:tauparameters} shows good agreement, as you can see in the last column of Table \ref{tab:gpuresults}.

The corresponding serial cost is simpler, given by
\begin{equation}
  C_C = n \taufcpu.
  \label{serialcostmodel}
\end{equation}
The speedup $ C_C / C_G $ can then be written
\begin{equation}
  \frac{C_C}{C_G} = \frac{n \taufcpu}{Mn \tau_F / N + N \tau_N + \tau_K} \approx \frac{Nn \kappa_F}{Mn + N^2 \kappa_N}. 
\end{equation}
where we have ignored the lower order term $ \tau_K $ because we are interested in the extreme case of large $ n $ and $ N $ and we have defined $ \kappa_F = \taufcpu / \tau_F $ and $ \kappa_N = \tau_N / \tau_F $.  In order to balance time discretization with interpolation errors, we will assume $ M \approx n^\alpha $ for $ \alpha < 1 $.  In fact in model problem with Chebyshev interpolation we could choose $ M \approx \log(n) $ but we will not assume we have such an exponentially convergent scheme in general.  Then if for fixed problem size $ n $ we choose $ N^2 \approx \kappa_N^{-1} n^{\alpha + 1} $ we get
\[
  \frac{C_C}{C_G} \approx \frac{\kappa_N^{-1/2} n^{(3 + \alpha)/2} \kappa_F}{n^{\alpha+1} + n^{\alpha + 1}} = \left( \frac{ \kappa_F} {2 \kappa_N^{1/2}} \right) n^{(1-\alpha)/2}
\]
which implies that for $ \alpha < 1 $, the larger the problem is, the more speedup we can get.

Similar to the cost model analysis of a parareal--type method in \cite[Section 3]{farhattime}, we see that the traditional parallel efficiency for our method is not as good as usual space--based parallelism.  However, the goal is to make the best use of an available GPU that otherwise would not be doing useful work, and we see that it is possible to get speedup in this case.

\section{Error analysis}
\label{analysis}

Nievergelt proves some error bounds on his approach in \cite{nievergeltparallel}, using forward Euler and linear interpolation.  We want to say something about the error in a slightly more abstract setting, with correspondingly abstract assumptions.

To begin, fix $ j $ and consider the interval $ \tdomain_j = [T_{j-1}, T_j] $.  Let $ g^j $ be the exact solution operator that maps an initial condition at $ t = T_{j-1} $ to a final condition at $ t = T_j $, and let $ g^j_{\Delta t} $ be a standard time-discretization method that approximates $ g^j $.  We will assume that our problem is well posed in the sense that $ g^j $ is Lipschitz as a function of the initial conditions,
\begin{equation}
  | g^j(\xi_1) - g^j(\xi_2) | \leq K | \xi_1 - \xi_2 |, \quad \xi_1, \xi_2 \in \Xi.
  \label{flipschitz}
\end{equation}

We evaluate $ g^j_{\Delta t} $ for many initial conditions in order to construct an interpolating function $ \phi^j_{\Delta t} $ which approximates the map $ g^j_{\Delta t} $ using some standard interpolation technique.  (The map $ g^j_{\Delta t} $ in turn approximates the exact operator $ g^j $, which is of course what we really want.)  Denote the true solution at $ t = T_{j-1} $ by $ \lambda_{j-1}^* $, but of course we will instead have an approximation to this value communicated from a previous processor, denoted by $ \lambda_{j-1} $.  We will also use $ g $ and $ g_{\Delta t} $ without superscripts to denote propagators on the whole time interval $ \tdomain = [T_0, T_N] $, and $ \phi_{\Delta t} $ will be the Nievergelt method used over the entire time domain (including several interpolations).

We begin our analysis by considering the error at the end of the fixed time domain $ \tdomain_j $:
\begin{align}
  | \phi^j_{\Delta t} (\lambda_{j-1}) - g^j(\lambda_{j-1}^*) | \leq | \phi^j_{\Delta t} (\lambda_{j-1}) &- g^j_{\Delta t}(\lambda_{j-1}) | + | g^j_{\Delta t}(\lambda_{j-1}) - g^j(\lambda_{j-1}) | \notag \\
  & + | g^j(\lambda_{j-1}) - g^j(\lambda_{j-1}^*) |,
  \label{errorsplitting}
\end{align}
where the first term is recognized as a standard interpolation error, the second term is a standard time discretization error, and the third term depends on how accurately the solution was computed on the previous time-slice $ \tdomain_{j-1} $.  Beginning with this last term and using \eqref{flipschitz}, we have
\[
  |g^j(\lambda_{j-1}) - g^j(\lambda_{j-1}^*) | \leq K | \lambda_{j-1} - \lambda_{j-1}^* | = K | \phi_{\Delta t}^{j-1} (\lambda_{j-1}) - g^{j-1} (\lambda_{j-2}^*) |,
\]
which contains exactly the same type of expression as \eqref{errorsplitting}, so that we can use it recursively to get an expression for the final error,
\begin{align*}
  | \phi^N_{\Delta t} (\lambda_{N-1}) - g^N(\lambda_{N-1}^*) | &\leq NK \left( \max_j | \phi^j_{\Delta t} (\lambda_{j-1}) - g^j_{\Delta t}(\lambda_{j-1}) | \right. \notag \\
  &\qquad + \left. \max_j | g^j_{\Delta t}(\lambda_{j-1}) - g^j(\lambda_{j-1}) | \right).
\end{align*}

We will require that our time discretization method satisfies
\begin{equation}
  | g_{\Delta t}^j(\lambda_{j-1}) - g^j(\lambda_{j-1}) | \leq C_1 \frac{1}{N} O(\Delta t^k),
  \label{timecondition}
\end{equation}
where $ C_1 $ is independent of the time interval $ T_j - T_{j-1} $ and $ k $ is the order of error in the underlying time discretization method,
\[
  | g_{\Delta t}(\lambda_0^*) - g(\lambda_0^*) | \leq D_1(T) O(\Delta t^k).
\]
The standard error bounds for many time discretization methods (including linear multistep methods) have $ D_1 $ growing {\em exponentially} with the length of the time interval that we are integrating over, so that \eqref{timecondition} is not hard to satisfy \cite{atkinsonnumerical,henricidiscrete}.

The analogous condition for the interpolation is easy to state but perhaps harder to interpret.  It is
\begin{equation}
  | \phi_{\Delta t}^j (\lambda_{j-1}) - g_{\Delta t}^j (\lambda_{j-1}) | \leq C_2 \frac{1}{N} O(\Delta \xi^\ell),
  \label{interpolationcondition}
\end{equation}
where $ \ell $ again represents the order of error of the underlying interpolation on the initial value space:
\[
  | \phi_{\Delta t} (\lambda_0^*) - g_{\Delta t} (\lambda_0^*) | \leq D_2(T) O(\Delta \xi^\ell).
\]
Here we are assuming that the interpolating function to the time discretization operator is more accurate if the time integration is over a short interval than if it is over a long interval (cf. \cite{nievergeltparallel}).  If we use polynomial interpolation, standard results give us 
\[
  | \phi_{\Delta t}^j (\lambda_{j-1}) - g_{\Delta t}^j (\lambda_{j-1}) | \leq C_3 \Delta \xi^{p+1} \max_{\xi \in \Xi} \left| \frac{\partial^{p+1} g_{\Delta t}^j(\xi)}{\partial \xi^{p+1}} \right|,
\]
where $ p $ is the order of the polynomial.  The quantity in absolute value can be interpreted as the $ D_2(T) $ above, and in general can be expected to grow very quickly in $ T $, so that \eqref{interpolationcondition} is also a very plausible assumption for a wide range of problems---in particular it certainly holds for our simple model problem \eqref{scalarmodel} with linear interpolation.

We put these simple results together in the following
\begin{theorem}
  Under the assumptions \eqref{flipschitz}, \eqref{timecondition} and \eqref{interpolationcondition}, we have
  \begin{equation}
    | \phi^j_{\Delta t} (\lambda_{j-1}) - g^j(\lambda_{j-1}^*) | \leq O(\Delta \xi^\ell) + O(\Delta t^k),
  \end{equation}
  and in particular if we balance interpolation and time discretization errors, the Nievergelt approach recovers the error of the underlying time discretization.
\end{theorem}

For linear differential equations, the solution is also linear as a function of its initial condition, so that the interpolation error can be made zero using only enough points to construct a basis for the initial condition space, as shown in \cite{nievergeltparallel}.  In the remainder of the paper we will focus on this simpler case.

We close this section with a few comments on the cost required to achieve this error.  The Nievergelt algorithm is very computationally costly in conventional terms.  In a sense, we trade off everything else against communication costs, assuming those to be dominant.  The result is that floating point operation counts and storage costs are very high.  For example, when applied to a nonlinear PDE discretized with $ r $ unknowns, the initial value space will have dimension $ r $.  If you sample uniformly in this high--dimensional space, with $ \ell $ points for each of the $ r $ unknowns, then the total number of operations to construct the propagation operator (or storage to store the final conditions) will grow as $ \ell^r $.  If parallel operations are truly free we may be willing to pay this cost, or in some cases we may be able to sample much more efficiently than uniformly.

In addition to the computational cost, in many dimensions the choice of how to sample the initial value space and how to perform the interpolation step is not straightforward.  If some suitable interpolation can be applied, than the Nievergelt method can in principle be applied and will result in parallel speedup, but in many cases this will not be practical and some other algorithm must be sought.

However, in the remainder of this paper we will consider linear problems, where the storage and cost issues can be greatly reduced and the interpolation issue disappears completely.  In particular, if after space discretization we have the system
\begin{equation}
  \frac{d}{dt} \bfy = A(t) \bfy + \bfb(t),
  \label{lineardiscrete}
\end{equation}
then the solution is an affine function of the initial condition $ \bfy_0 $.  With $ M $ points in the spatial discretization, then $ \bfy_0 = \sum_{j=1}^M y_j e_j $ for the unit basis vectors $ e_j $, and we can solve \eqref{lineardiscrete} once with a zero initial condition and solve a homogeneous version of \eqref{lineardiscrete} $ M $ times, once for each basis vector.  This can be done with no knowledge of the true initial condition coefficients $ y_j $.  Then we reconstruct the solution
\begin{equation}
  \bfy = g_{\Delta t}(0) + \sum_{j=1}^M y_j g_{\Delta t}(e_j).
  \label{linearreconstruct}
\end{equation}
We use the notation $ g_{\Delta t} $ (rather than $ \phi_{\Delta t} $) because there is no need for interpolation in this linear setting and this technique introduces no additional error beyond rounding error.

\section{The heat equation and wave equation}
\label{pdesection}

Having examined the Nievergelt approach numerically on a simple scalar initial value problem and done some analysis, we now turn to some results for slightly more complicated problems.  First we consider the variable coefficient heat equation
\begin{equation}
  \frac{\partial y}{\partial t} = \left( 1 + \frac{1}{4} \sin(t) \right) \frac{\partial^2 y}{\partial x^2} + b(t), \quad x \in [0,1]
  \label{heatequation}
\end{equation}
in one dimension with homogeneous Dirichlet boundary conditions and with $ b(t) $ and the initial conditions chosen so that the true solution is given by
\[
  y = \cos(t) \sin(\pi x).
\]
We note that \eqref{heatequation} has a nonconstant coefficient and so that methods on exponential quadrature or Laplace transforms \cite{hochbruckexponentialrk,gavrilyukexponentially,sheenparallel} and the paraexp method \cite{ganderparaexp} are not applicable.  

\begin{table}
  \tbl{Time to solution in seconds for the model problem \eqref{heatequation} using various numbers of processors $ N $.  Here $ \Delta t = 0.005 $ and $ \Delta x = 0.1 $, we are propagating from $ t = 0 $ to $ t = 10 $ using backward Euler for timestepping.}
  {\begin{tabular}{l|ll}
  $ N $& \ttotal& \tcomm \\
  \hline
  1&  0.2870& --- \\
  2& 1.0547& 1.71e-4 \\
  4& 0.5264& 4.13e-4 \\
  8& 0.2637& 4.51e-4 \\
  16& 0.1322& 2.38e-4 \\
  32& 0.06730& 6.16e-4 \\
  64& 0.03453& 6.71e-4 \\
 128& 0.01763& 6.13e-4 \\
  \end{tabular}}
  \label{heatlong}
\end{table}

After discretization in space using standard centered finite differences and $ M $ spatial points, \eqref{heatequation} can be written in the form \eqref{lineardiscrete}.  On each parallel time slice we solve this problem $ M + 1 $ times with different initial conditions (and with $ b(t) \neq 0 $ only once) and reconstruct the solution using \eqref{linearreconstruct} and a single communication step.

\begin{table}
  \tbl{Parareal comparison to Table \ref{heatlong}, with $ \Delta t = 0.005, \Delta T = 0.1 $, and propagating to a final time of 10.}
    {\begin{tabular}{l|ll|ll|ll|ll}
      & $k=2$& & $k=3$& & $k=4$& & $k=8$& \\
   $N$& \ttotal& \tcomm& \ttotal& \tcomm& \ttotal& \tcomm& \ttotal& \tcomm \\
      \hline
      2& 0.3047& 2.79e-5&   0.4533& 4.20e-5&   0.6028& 4.89e-5&   1.204& 1.06e-4  \\
      4& 0.1500& 2.88e-4&   0.2250& 5.48e-4&   0.2994& 4.72e-4&   0.5986& 8.81e-4  \\
      8& 0.07744& 6.22e-4&  0.1152& 7.53e-4&   0.1533& 6.56e-4&   0.3054& 1.03e-3  \\
     16& 0.04020& 1.35e-3&  0.05938& 1.87e-3&  0.07721& 5.66e-4&  0.1543& 7.54e-4  \\     
     32& 0.02163& 1.72e-3&  0.03119& 1.57e-3&  0.04093& 1.85e-3&  0.07940& 1.56e-3  \\
    \end{tabular}}
  \label{paraheat}
\end{table}

\begin{table}
  \tbl{Time to solution in seconds for the model problem \eqref{heatequation} using various numbers of processors $ N $.  Here $ \Delta t = 0.001 $ and $ \Delta x = 0.05 $, we are propagating from $ t = 0 $ to $ t = 10 $ using backward Euler for timestepping.}
{\begin{tabular}{l|ll}
  $ N $& \ttotal& \tcomm \\
  \hline
  1& 1.641& --- \\
  2& 12.309& 4.84e-4 \\
  4& 6.112& 2.26-3 \\
  8& 3.068& 1.11e-2 \\
  16& 1.530& 4.92e-3 \\
  32& 0.7748& 1.12e-2 \\
  64& 0.3859& 1.57e-3 \\
  128& 0.1939& 6.07e-4 \\
  256& 0.0984& 4.37e-4 \\
\end{tabular}}
\label{heatlongfine}
\end{table}

Some numerical results are shown in Table \ref{heatlong}, where we can see a speedup of about 15.  Note that we are able to get this speedup partially because we are willing to accept poor spatial resolution.  If more points in space are required, we will need many more processors to see speedup, but there is no real barrier to doing this provided the processors are available.  Scaling results for a parareal implementation of this problem are shown in Table \ref{paraheat}.  See also Tables \ref{heatlongfine} and \ref{heatlongfineqb} which achieve similar speedups for larger problems with more processors. 

\begin{table}
  \tbl{Time to solution in seconds for the model problem \eqref{heatequation} using various numbers of processors $ N $.  Here $ \Delta t = 0.0005 $ and $ \Delta x = 0.01 $, propagating from $ t = 0 $ to $ t = 10 $, using backward Euler for timestepping.}
  {\begin{tabular}{l|ll}
  $ N $& \ttotal& \tcomm \\
  \hline
  1& 16.166& --- \\
  32& 45.593& 0.230 \\
  64& 22.901& 0.140 \\
  128& 11.867& 0.203 \\
  256& 6.130& 0.361 \\
  512& 3.118& 0.0963 \\
  1024& 1.595& 0.0780 \\
  \end{tabular}}
  \label{heatlongfineqb} 
\end{table}

\begin{table}
  \tbl{Timing comparisons between parareal and the Nievergelt minimal communication method.  Here $ \Delta t = 0.005 $ and we use backward Euler and include an artificial communication delay of 0.001 seconds for every receive.  The parareal method uses $ \Delta T = 0.1 $.}
    {\begin{tabular}{l|ll|ll|ll|ll}
      & Nievergelt& & $k=2$& & $k=4$& & $k=8$& \\
   $N$& \ttotal& \tcomm& \ttotal& \tcomm& \ttotal& \tcomm& \ttotal& \tcomm\\
      \hline
      1& 0.2870&  --- \\
      2& 1.0690& 4.31e-3&  0.3251&  2.37e-2&  0.6349&  3.25e-2&   1.4513& 7.98e-2  \\
      4& 0.5360& 7.82e-3&  0.1680&  1.76e-2&  0.3568&  5.66e-2&   0.7220& 1.21e-1  \\
      8& 0.2834& 1.89e-2&  0.1105&  3.38e-2&  0.2049&  5.41e-2&   0.4199& 1.16e-1  \\
     16& 0.1500& 1.73e-2&  0.07022& 3.18e-2&  0.1347&  5.85e-2&   0.2512& 9.97e-2  \\     
     32& 0.08318& 1.64e-2& 0.05538& 3.55e-2&  0.1067&  6.75e-2&   0.2002& 1.23e-1  \\
     64& 0.04521& 1.15e-2& 0.04460& 3.44e-2&  0.08416& 6.43e-2&   0.1645& 1.25e-1  \\
    \end{tabular}}
  \label{heatdelay}
\end{table}

Comparisons of the communication times in Tables \ref{heatlong} and \ref{paraheat} show that we do not, in fact, live in a world where parallel floating point operations are free compared to communication time, and the Nievergelt algorithm is competitive with parareal only if we require an unrealistically large number of parareal iterations.  In Table \ref{heatdelay} we simulate additional latency by enforcing a delay of one millisecond for every message that is received from another processor, and it is clear that in this simulated high-latency environment the additional communication steps of the parareal method are a drawback, especially for problems where a large number of iterations might be required.  As parallel computations become cheaper relative to communication, the Nievergelt approach becomes more and more promising in comparison to the parareal method.

The parareal algorithm is reported to work relatively poorly for second--order problems and for hyperbolic problems \cite{farhattime,ganderanalysismodified,gandercharacteristics,ganderkrylov}, and many of the other time-parallel methods we considered in the introduction are restricted to parabolic problems \cite{hochbruckexponentialrk,gallopoulosefficient,sheenparallel}.  As our final numerical example we consider the Nievergelt method applied to these types of problems.  In \cite{ganderanalysismodified,ganderkrylov}, a version of the parareal algorithm for second--order and hyperbolic algorithms is presented---however, this method requires the communication of more initial conditions from more distant processors (not just the immediately preceding time slice), and so is not competitive in trying to minimize communication.

We will use an example problem from \cite{trefethenspectral}, Program 19, which is the wave equation 
\begin{equation}
  \frac{\partial^2 y}{\partial t^2} = \frac{\partial^2 y}{\partial x^2}, \quad x \in [0,1]
  \label{waveequation}
\end{equation}
with homogeneous Dirichlet boundary conditions and with the initial conditions such that the solution is a single peaked wave propagating to the left.  This problem is discretized in space with a spectral method on Chebyshev points, and time-stepping is with the explicit second-order leapfrog method.  We will use a varying number $ M $ of spatial points but always choose a timestep of $ \Delta t = 8 / M^2 $ to satisfy a CFL condition.  See \cite{trefethenspectral} for details for this model problem.

\begin{table}
  \tbl{Time to solution for model problem \eqref{waveequation} with $ M = 40 $ Chebyshev points, integrating to a final time of 16.}
    {\begin{tabular}{l|ll}
      $N$& \ttotal& \tcomm \\
      \hline
      1& 3.193& --- \\
      2& 64.950& 0.0782 \\
      4& 32.453& 0.00122 \\
      8& 16.334& 0.0907 \\
      16& 8.163& 0.0416 \\
      32& 4.124& 0.00121 \\
      64& 2.061& 0.0299 \\
      64& 2.049& 0.0288 \\
      128& 1.039& 0.0223 \\
      256& 0.552& 0.0620
    \end{tabular}}
  \label{wave40}
\end{table}

The time-domain decomposition proceeds much as it did for the heat equation, the only important difference being that we have two initial conditions (data from the previous timestep and the step before that) to deal with and to transfer from subdomain to subdomain.  This has the effect of doubling the initial value space, that is, if we have $ M = 80 $ spatial points, we need to propagate 160 initial conditions to get the full map from initial conditions to final conditions.

Numerical results are shown in Tables \ref{wave40} and \ref{wave80}.

\begin{table}
  \tbl{Time to solution for model problem \eqref{waveequation} with $ M = 80 $ Chebyshev points, integrating to a final time of 20.}
    {\begin{tabular}{l|ll}
      $N$& \ttotal& \tcomm \\
      \hline
      1& 9.058& --- \\
      16& 90.969& 0.399 \\
      32& 45.673& 0.291 \\
      64& 22.845& 0.150 \\
      128& 11.454& 0.104 \\
      256& 5.809& 0.211 \\
      512& 3.030& 0.212 \\
    \end{tabular}}
  \label{wave80}
\end{table}

\section{Conclusion}
\label{conclusion}

We have examined the method of Nievergelt for the parallel solution of time-dependent problems, a method which essentially builds the time propagation operator on each time-slice in parallel with no knowledge of the initial conditions.

It is clear from a viewpoint of counting arithmetic operations that the Nievergelt approach is not efficient, and it may never be suitable for nonlinear problems with a large number of degrees of freedom.  However, we are moving to a future where computational cores are essentially free---hundreds of them may be sitting unused on your desktop, and if there is any way to use them profitably it may be worthwhile even if the efficiency is poor.  In addition, for some applications the total time to solution is the only concern, and the amount of computational resources spent is less important.  It is this kind of hardware environment and these kind of applications that the parareal algorithm was invented to deal with.  The Nievergelt method compares well to the parareal algorithm but has the important advantage of always requiring less communication, and it is applicable to more situations than methods based on exponential quadrature or Laplace transforms.

In our work above, we have shown speedup of 3 on 32 processors for a scalar ODE problem, a speedup of 10 on 1024 processors for a variable coefficient heat equation and speedup of 5 on 256 processors for a wave equation.  These results are somewhat worse than those in published parareal implementations, but roughly comparable.  The original parareal paper \cite{lionsresolution} predicted speedups of 8 to 18, and in \cite{liumodified} the authors get speedup of 5 or 6 on 200 processors.  The best speedups in \cite{farhattime} and \cite{farhatimplicit} are from 4 to 6.  For the different paraexp method, speedups are reported to be from 4 to 6 \cite{ganderparaexp}.  In addition to speedup on traditional parallel clusters, we have seen that implementing the Nievergelt algorithm on a GPU provides modest speedup, with this case being of particular interest because many computational scientists will have a GPU sitting idle in their desktop while they wait for a CPU to compute their results.

The assumptions that we are using in this work, namely that parallel computation is so cheap as to be nearly free when compared to communication costs, are not satisfied on traditional parallel machines.  However, to some extent this assumption holds on modern graphics hardware, and we have shown that algorithms can be designed for a future world where this kind of pattern is more common.  Even if the very high computational and storage costs of the Nievergelt method make it impractical in many cases, we argue that something like it could influence at least part of the design of future parallel methods and that understanding this extreme case may point the way to exploiting future hardware where communication costs dominate computational costs.  Aside from the model problem at the beginning, we have focused on linear problems, because in this case it is easier to construct and represent the map from initial conditions to final conditions.  In the general nonlinear case, using this map will require approximating it and using some high-dimensional interpolation scheme, which may be prohibitively expensive.  However, there is no additional communication necessary for nonlinear problems.

We have divided our time domain into slices and assigned each slice to a single processor on a parallel machine, but it is worth noting that there is a great deal more parallelism inherent in the problem than we have exploited.  The map from initial to final conditions is constructed by integrating in time a large number of initial conditions given by some discretization of the initial condition space, and in principle each of these time integrations could be done on a separate processor or processing core in parallel.  The lack of need for synchronization in the whole algorithm (except at the final step) suggests the possibility of a distributed computing model for some applications.

There are several possible directions for future research.  One is to combine the Nievergelt approach with a parareal method in some way, to get a balance of the communication costs versus the greatly added computation of Nievergelt.  Another is to do more realistic nonlinear problems so as to more carefully consider the role of interpolation in the method.  Ongoing research is to combine the traditional parallel implementation with the GPU implementation to see the potential of the method on heterogeneous systems.

\section*{Acknowledgements}

This work was supported in part by the U.S. National Science Foundation under Grant Number DMS-07-39382.  Computational resources were provided by the Louisiana Optical Network Initiative and the Stellar Group at the Center for Computational Technology at Louisiana State University.

\bibliographystyle{gCOM}
\bibliography{minimal}

\begin{thebibliography}{10}
\providecommand{\url}[1]{\texttt{#1}}
\providecommand{\urlprefix}{URL }

\bibitem{atkinsonnumerical}
K.E. Atkinson, \emph{An introduction to numerical analysis}, Wiley (1989).

\bibitem{christliebparallel}
A.J. Christlieb, C.D. MacDonald,  and B.W. Ong, \emph{Parallel high-order
  integrators}, SIAM J. Sci. Comput. 32 (2010), pp. 818--835.

\bibitem{duarteparareal}
M. Duarte, M. Massot,  and S. Descombes, \emph{Parareal operator splitting
  techniques for multi-scale reaction waves: numerical analysis and
  strategies}, Math. Model. Numer. Anal. 45 (2011), pp. 825--852.

\bibitem{farhattime}
C. Farhat and M. Chandesris, \emph{Time-decomposed parallel time-integrators:
  theory and feasibility studies for fluid, structure, and fluid--structure
  applications}, Int. J. Numer. Meth. Engng. 58 (2003), pp. 1397--1434.

\bibitem{farhatimplicit}
C. Farhat, J. Cortial, C. Dastillung,  and H. Bavestrello, \emph{Time-parallel
  implicit integrators for the near-real-time prediction of linear structural
  dynamic responses}, Int. J. Numer. Meth. Engng. 67 (2006), pp. 697--724.

\bibitem{fischerparareal}
P.F. Fischer, F. Hecht,  and Y. Maday, \emph{A parareal in time semi-implicit
  approximation of the {N}avier-{S}tokes equations}, in \emph{Domain
  Decomposition Methods in Science and Enginnering}, e.a. Timothy J.~Barth,
  ed., \emph{Lecture Notes in Computational Science and Engineering}, vol.~40,
  Springer,  2005, pp. 443--440.

\bibitem{gallopoulosefficient}
E. Gallopoulos and Y. Saad, \emph{Efficient solution of parabolic equations by
  {K}rylov approximation methods}, SIAM J. Sci. Stat. Comput. 13 (1992), pp.
  1236--1264.

\bibitem{ganderanalysismodified}
M. Gander and M. Petcu, \emph{Analysis of a modified parareal algorithm for
  second-order ordinary differential equations}, AIP Conference Proceedings 936
  (2007), pp. 233--236.

\bibitem{ganderkrylov}
---{}---{}---, \emph{Analysis of a {K}rylov subspace enhanced parareal
  algorithm for linear problems}, ESAIM Proc. 25 (2008), pp. 114--129.

\bibitem{gandercharacteristics}
M.J. Gander, \emph{Analysis of the parareal algorithm applied to hyperbolic
  problems using characteristics}, Bol. Soc. Esp. Mat. Apl. 42 (2008), pp.
  5--19.

\bibitem{ganderparaexp}
M.J. Gander and S. G\"uttel, \emph{Paraexp: {A} parallel integrator for linear
  initial-value problems}, SIAM J. Sci. Comput.  (2013), p. to appear.

\bibitem{ganderanalysis}
M.J. Gander and S. Vandewalle, \emph{Analysis of the parareal time-parallel
  time-integration method}, SIAM J. Sci. Comput. 29 (2007), pp. 556--678.

\bibitem{gavrilyukexponentially}
I.P. Gavrilyuk and L. Makarov, \emph{Exponentially convergent algorithms for
  the operator exponential with applications to inhomogeneous problems in
  {B}anach spaces}, SIAM J. Numer. Anal. 43 (2005), pp. 2144--2177.

\bibitem{henricidiscrete}
P. Henrici, \emph{Discrete variable methods in ordinary differential
  equations}, Wiley (1961).

\bibitem{hochbruckexponentialrk}
M. Hochbruck and A. Ostermann, \emph{Exponential {R}unge-{K}utta methods for
  parabolic problems}, Appl. Numer. Math. 53 (2005), pp. 323--339.

\bibitem{lionsresolution}
J.L. Lions, Y. Maday,  and G. Turinici, \emph{R\'esolution d'{EDP} par un
  sch\'ema en temps parar\'eel}, C.R. Acad. Sci. Paris 332 (2001), pp.
  661--668.

\bibitem{liumodified}
Y. Liu and J. Hu, \emph{Modified propagators of parareal in time algorithm and
  application to {P}rinceton {O}cean model}, Int. J. Numer. Meth. Fluids 57
  (2008), pp. 1793--1804.

\bibitem{madayparareal}
Y. Maday and G. Turinici, \emph{The parareal in time iterative solver: a
  further direction to parallel implementation}, in \emph{Domain Decomposition
  Methods in Science and Engineering}, e.a. Timothy J.~Barth, ed.,
  \emph{Lecture Notes in Computational Science and Engineering}, vol.~40,
  Springer,  2005, pp. 441--448.

\bibitem{minionhybrid}
M.L. Minion, \emph{A hybrid parareal spectral deferred corrections method},
  Comm. App. Math. and Comp. Sci. 5 (2010), pp. 265--301.

\bibitem{nievergeltparallel}
J. Nievergelt, \emph{Parallel methods for integrating ordinary differential
  equations}, Comm. ACM 7 (1964), pp. 731--733.

\bibitem{sheenparallel}
D. Sheen, I.H. Sloan,  and V. Thom\'ee, \emph{A parallel method for time
  discretization of parabolic equations based on {L}aplace transformation and
  quadrature}, IMA J. Numer. Anal. 23 (2003), pp. 269--299.

\bibitem{trefethenspectral}
L.N. Trefethen, \emph{Spectral Methods in {MATLAB}}, SIAM (2000).

\end{thebibliography}

\end{document}